\documentclass{amsart}

\usepackage{latexsym,amsmath,amssymb,amsfonts,amscd,graphics,appendix,amsxtra}
\usepackage[mathscr]{eucal}
\usepackage[all]{xypic}
\usepackage[normalem]{ulem}
\usepackage[pdftex]{graphicx}

\numberwithin{equation}{section}
\theoremstyle{plain}
\newtheorem{theorem}[equation]{Theorem}

\newtheorem{proposition}[equation]{Proposition}
\newtheorem{corollary}[equation]{Corollary}

\theoremstyle{remark}
\newtheorem{remark}[equation]{Remark}

\theoremstyle{definition}

\newtheorem{example}[equation]{Example}

\def\Spec{\operatorname{Spec}}

\def\Hom{\operatorname{Hom}}

\def\Sym{\operatorname{Sym}}

\newcommand{\bP}{\mathbb{P}}
\newcommand{\bA}{\mathbb{A}}
\newcommand{\bG}{\mathbf{G}}
\newcommand{\bB}{\mathbf{B}}
\newcommand{\bK}{\mathbf{K}}

\newcommand{\bD}{\mathbf{D}}
\newcommand{\bZ}{\mathbb{Z}}
\newcommand{\bQ}{\mathbb{Q}}
\newcommand{\bR}{\mathbb{R}}
\newcommand{\bC}{\mathbb{C}}

\newcommand{\fX}{\mathfrak{X}}

\newcommand{\calA}{\mathcal{A}}

\newcommand{\Aut}{\mathrm{Aut}}

\newcommand{\Hilb}{\mathrm{Hilb}}

\newcommand{\Proj}{\mathrm{Proj}}

\newcommand{\Gr}{\mathrm{Gr}}

\newcommand{\gquot}{/\!\!/}

\newcommand{\SO}{\mathrm{SO}}
\newcommand{\Sp}{\mathrm{Sp}}

\newcommand{\SL}{\mathrm{SL}}

\newcommand{\calB}{\mathcal{B}} 
\newcommand{\calD}{\mathcal{D}} 
\newcommand{\calH}{\mathcal{H}} 
 
\newcommand{\calO}{\mathcal{O}}
\newcommand{\calL}{\mathcal{L}}

\newcommand{\calF}{\mathcal{F}}
\newcommand{\calh}{\mathfrak{h}}

\newcommand{\calM}{\mathcal{M}}
\newcommand{\calP}{\mathcal{P}}
\newcommand{\calN}{\mathcal{N}} 
 
\newcommand{\calC}{\mathcal{C}}
\newcommand{\calX}{\mathcal{X}}

\newcommand{\calV}{\mathscr{V}}

\newcommand{\bdg}{\mathbb D/\Gamma}

\begin{document}
\bibliographystyle{amsalpha}
\title[Perspectives on Moduli Spaces]{Perspectives on the construction and compactification of moduli spaces}

\author{Radu Laza}
\address{Stony Brook University, Department of Mathematics, Stony Brook, NY 11794}
\email{rlaza@math.sunysb.edu}
\thanks{The author was partially supported by NSF grants DMS-1200875 and DMS-125481 (CAREER), and by a Sloan Fellowship.}

\begin{abstract}
 In these notes, we introduce various approaches (GIT, Hodge theory, and KSBA) to constructing and compactifying moduli spaces. We then discuss the pros and cons for each approach, as well as some connections between them. 
\end{abstract}

\maketitle

\section*{Introduction}
A central theme in algebraic geometry is the construction of compact moduli spaces with geometric meaning. The two early successes of the moduli theory - the construction and compactification of the moduli spaces of curves $\overline M_g$ and principally polarized abelian varieties (ppavs) $\overline \calA_g$ - are models that we try to emulate. While very few other examples are so well understood, the tools developed to study other moduli spaces have led to new developments and unexpected directions in algebraic geometry. The purpose of these notes is to review three standard approaches to constructing and compactifying moduli spaces: GIT, Hodge theory, and MMP, and to discuss various connections between them. 

One of the oldest approach to moduli problems is Geometric Invariant Theory (GIT). The idea is natural: The varieties in a given class can be typically embedded into a fixed projective space. Due to the existence of the Hilbert schemes, one obtains a quasi-projective variety $X$ parametrizing embedded varieties of a certain class. Forgetting the embedding amounts to considering the quotient $X/G$ for a certain reductive algebraic group $G$. Ideally, $X/G$ would be the moduli space of varieties of the given class. Unfortunately, the naive quotient $X/G$ does not make sense; it has to be replaced by the GIT quotient $X\gquot G$ of Mumford \cite{GIT}. While $X\gquot G$ is the correct quotient from an abstract point of view, there is a price to pay: it is typically difficult to understand which are the semistable objects (i.e. the objects parameterized by $X\gquot G$) and then some of the semistable objects are too degenerate from a moduli point of view. Nonetheless, $X\gquot G$ gives a projective model for a moduli space with weak modular meaning. Since the GIT model $X\gquot G$ is sometimes more accessible than other models, $X\gquot G$ can be viewed a first approximation of more desirable compactifications of the moduli space. 

A different perspective on moduli is to consider the variation of the cohomology of the varieties in the given moduli stack $\calM$. From this point of view, one considers the homogeneous space $\mathbb D$  that classifies the Hodge structures of a certain type, and then the quotient  $\mathbb D/\Gamma$ which corresponds to forgetting the marking of the cohomology. The ideal situation would be a period map $\calP:\calM\to \mathbb D/\Gamma$ which is an isomorphism, or at least a birational map. Results establishing the  (generic) injectivity of the period map are called ``Torelli theorems'', and a  fair number of such results are known. Unfortunately,   the image of $\calP$ in $\mathbb D$ is typically very hard to understand: Griffiths' transversality says that the periods of algebraic varieties vary in a constrained way, which gives a highly non-trivial systems of differential equations. Essentially, the only cases where we don't have to deal with these issues are the classical cases of ppavs and $K3$ surfaces, for which all our knowledge on their moduli is obtained by this Hodge theoretic construction. Furthermore, having a good period map gives numerous geometric consequences. The reason for this is that the spaces $\bdg$ have a lot of structure that can be translated into geometric properties.
 While it is advantageous  to get a description of the moduli space as a locally symmetric variety, in practice very few examples are known. We will briefly mention some enlargement of the applicability of period map constructions to moduli beyond ppavs and $K3$s. Finally, we will review some work of Looijenga which gives some comparison results for the case when both the GIT and Hodge theoretic approach are applicable. This is in some sense an ideal situation as both geometric and structural results exist. 

While the first two approaches are based on considering the properties of smooth objects, and then constructing a global moduli space. The third approach takes a different tack:  one constructs a moduli space by gluing local patches. This gives a moduli stack, and the main issue is to carefully choose degenerations such that one obtains a proper and separated stack. By the valuative criteria, it suffices to consider $1$-parameter degenerations. From a topological point of view,  the ideal model is  a semistable degeneration $\fX/\Delta$, but then the central fiber is far from unique. A fundamental insight comes from the minimal model program (MMP): the canonical model of varieties of general type is unique. Consequently, by allowing ``mild'' singularities, one obtains a unique limit for a $1$-parameter degeneration, leading to a proper and separated moduli stack, and (under mild assumptions) even a projective coarse moduli space. This theory was developed by Koll\'ar--Shepherd-Barron--Alexeev (KSBA) with contributions from other authors. The relationship between the KSBA approach and the other two approaches is not well understood. We briefly review some partial results on this subject. In one direction, the connection between du Bois and semi-log-canonical (slc) singularities gives a link between KSBA and Hodge theory.  In a different direction, the Donaldson--Tian theory of K-stability establishes a connection between GIT and KSBA stability. These topics are rapidly evolving and suggest that much is still to be explored in moduli theory. 

The overarching theme of these notes is that each approach sheds light on a different aspect of the moduli problem under consideration. By taking together different approaches one obtains a fuller picture of a moduli space and its compactifications.  As examples of this principle, we point out to the theory of variation of GIT (VGIT) quotients of Thaddeus and Dolgachev--Hu, and the study of log canonical models for $\overline{M}_g$, the so called Hassett--Keel program. 

\subsection*{Disclaimer} These notes reflect the interests and point of view of the author. We have tried to give a panoramic view of a number of topics in moduli theory and to point out some relevant references for further details. In particular, we point out the related surveys \cite{lazavgit}, \cite{kollar}, \cite{yanosur}. We apologize for any omissions and inaccuracies. For instance, there is no discussion of deformation theory (\cite{sernesi}, \cite{defth}), stacks (\cite{dejong}), or the log geometry point of view (\cite{abramovich}). 

\subsection*{Acknowledgement} 
I have benefited from discussions with many people (including V. Alexeev, S. Casalaina-Martin, R. Friedman, P. Hacking, B. Hassett, E. Looijenga) over the years. I am particularly grateful to Y. Odaka and Z. Patakfalvi for some key comments on an earlier version of these notes. 
\section{The GIT Approach to constructing moduli spaces}\label{sectgit}
Geometric Invariant Theory (GIT) is probably the most natural and classical approach to constructing moduli spaces. In this section, we will review some of the main points of the GIT approach and survey some applications of GIT to moduli. The standard reference for GIT is Mumford et al. \cite{GIT}. Other good textbook references for GIT include: \cite{newstead}, \cite{dolgachev}, and \cite{mukaib}. For an extended discussion of the material included in this section see the survey \cite{lazavgit}. 

\subsection{Basic GIT and Moduli}
\subsubsection{Many moduli spaces are naturally realized as quotients $X/G$, where  $X$ is some (quasi-)projective variety and $G$ a reductive algebraic group acting on $X$.}\label{git1} The following results lead to presentations of moduli spaces as quotients $X/G$:

\medskip

\noindent (1) Given a class of polarized varieties, it is typically possible to give a uniform embedding result: {\it for all $(V,L)$ in the given class, for $k$ large and divisible enough (independent of $V$), $L^k$ is very ample and embeds $V$ into a fix projective space $\bP^{N}$ ($N=N(k)$ independent of $V$)}. For example:
\begin{itemize}
\item (Bombieri's Theorem) {\it For $V$ a smooth surface with $K_V$ big and nef, for all $k\ge 5$, the linear system $|k K_V|$ gives a birational morphism $V\to \overline V\subset \bP^N$, where $\overline V$ is the normal surface obtained by contracting all the smooth $(-2)$-curves orthogonal to $K_V$.} (Note, $\overline V$ is the canonical model of $V$ and has at worst du Val singularities.) 
\item (Mayer's Theorem) {\it Let $D$ be a big and nef divisor on a $K3$ surface, for all $k\ge 3$, the linear system $|k D|$ gives a birational morphism $V\to \overline V\subset \bP^N$, where $\overline V$ is the normal surface obtained by contracting all the smooth $(-2)$-curves orthogonal to $D$.}
\end{itemize}
While similar results for singular varieties and higher dimensions are more subtle, satisfactory (but non-effective) results exist in high generality (e.g. see \cite{alexeevbound}, \cite{haconxu}).  In other words, we can assume wlog that all varieties $V$ in a certain class are embedded in a fixed projective space $\bP^N$.  

\medskip

\noindent (2) There exists a {\it fine} moduli space for embedded schemes $V$ in $\bP^N$ with fixed numerical invariants (i.e. Hilbert polynomial $p_V(t)$): it is the Hilbert scheme $H:=\mathrm{Hilb}_{p_V}(\bP^N)$. This is a well known story (e.g. \cite[Ch. 1]{kollarrat}), we only want to emphasize here the connection between flatness and the preservation of numerical invariants (\cite[Thm. III.9.9]{har}), and the fact that the Hilbert scheme is one of the very few instances of fine moduli spaces in algebraic geometry. 

\medskip

\noindent (3) To pass from the Hilbert scheme $H$ to a moduli space, there are two final steps. First, the Hilbert scheme parameterizes many objects that might have no connection to the original moduli problem (e.g. entire components of $H$ might parameterize strange non-reduced schemes). Thus, we need to restrict to the locus $X\subset H$ of ``good'' objects. To get a good theory, it is needed that $X$ is a locally closed subvariety of $H$ (e.g. the defining conditions for good objects are either open or closed conditions). Local closeness for the moduli functor holds quite generally, but sometimes it is quite subtle (e.g. \cite{husks}). Finally, to forget the embedding $V\subset \bP^N$ amounts to allowing linear changes coordinates on $\bP^N$. In conclusion, by this construction, we essentially obtained a moduli space for varieties of a given class as a global quotient $X/G$, where $X\subset H$ is as above and $G=\mathrm{PGL}(N+1)$. (For technical reasons, we replace $X$ by its closure and $\mathrm{PGL}(N+1)$ by $\mathrm{SL}(N+1)$ in what follows.)

\subsubsection{The naive quotient $X/G$ typically doesn't make sense. The correct solution is the GIT quotient $X\gquot G$.}\label{git2} Ideally, we would like that 
\begin{itemize}
\item[(a)] {\it the quotient $X/G$ gives $1$-to-$1$ parameterization of the $G$-orbits in $X$}, and  
\item[(b)] {\it $X/G$ has the structure of an algebraic variety} (s.t. $X\to X/G$ is a morphism that is constant on orbits). 
\end{itemize}
This is rarely possible as the example of $G=\bC^*$ acting in the standard way on $\bA^1=\bC$ shows: There are two  orbits, $\bA^1\setminus\{0\}$ and $\{0\}$, but they cannot give two separate points in $X/G$ as this would contradict the continuity of $X\to X/G$ (N.B. $\{0\} \subset\overline{\bA^1\setminus\{0\}}=\bA^1$). The GIT solution is to relax the condition (a) and then use (b) to define a quotient in a universal categorical sense. In the affine case, $X=\Spec R$, it is easy to see that there is  only one possible choice 
$$X/G:=\Spec R^G,$$
where $R^G$ is a ring of $G$-invariant regular functions (automatically a finitely generated algebra if $G$ is reductive). With this definition, all the expected properties of $X\to X/G$ hold except for (a), which is replaced by \begin{itemize}
\item[(a')]{\it every point in $X/G$ corresponds to a unique closed orbit} (and two orbits map to the same point in $X/G$ iff the intersection of their closures is non-empty). 
\end{itemize}
For example, the invariant ring $\bC[x]^{\bC^*}$ (with $t\in \bC^*$ acting by $x\to tx$) is the constants $\bC$ and thus $\bA^1/\bC^*=\{\ast\}$ corresponding to two different orbits, one of which (i.e. $\{0\}$) is closed. 

In general, a quotient $X/G$ can be constructed by gluing quotients of open affine $G$-invariant neighborhoods of points in $X$. For simplicity, we restrict here to the case of $X$ being a projective variety with an ample $G$-linearized line bundle $\calL$. In this situation, the correct quotient from the point of algebraic geometry is 
\begin{equation}\label{defgit}
X\gquot G:=\Proj R(X,\calL)^G,
\end{equation}
where $R(X,\calL)=\oplus_n H^0(X,\calL^n)$. Note that in the projective situation, the natural map $X\dashrightarrow X\gquot G$ is only a rational map: it is defined only for {\it semistable points $x\in X^{ss}$}, i.e. points for which there exists an invariant section $\sigma\in H^0(X,\calL^n)$ non-vanishing at $x$. The stable locus $X^{s}\subseteq X^{ss}$ is the (open) set of semistable points $x$ for which the orbit $G\cdot x$ is closed in $X^{ss}$ and the stabilizer $G_x$ is finite. The quotient $X^s/G$ is  a geometric quotient, i.e. satisfies both conditions (a) and (b), and thus a good outcome for a moduli problem. 

At this point, we already see some issues with constructing a moduli space via GIT. First, the set of semistable points is somewhat mysterious and might not be what is expected. Secondly, if there exist strictly semistable points (i.e. $X^{ss}\setminus X^{s}\neq \emptyset$), then several orbits will correspond to the same point in the quotient $X\gquot G$. Thus, usually, the GIT quotients are not ``modular'' at the boundary.

\begin{remark}
It is well known that it is essential to work with reductive groups $G$, in order to obtain finitely generated rings of invariants $R^G$. However,  there exist natural situations when $G$ is not reductive. We point  to \cite{knonred} for some techniques to handle these cases. For some  concrete examples of non-standard GIT (e.g. $G$ non-reductive or non-ample linearization $\calL$) see \cite{g4ball,g4git}.
\end{remark}

\subsubsection{The GIT quotient $X\gquot G$ depends on the choice of linearization. This gives flexibility to the GIT construction, which is sometimes very useful.}\label{git3} By definition (cf. \eqref{defgit}), the GIT quotient depends on the choice of linearization $\calL$, so it is more appropriate to write $X\gquot_\calL G$. A surprising fact discovered by Dolgachev-Hu \cite{dh} and Thaddeus \cite{thaddeus} is that the dependence on $\calL$ is very well behaved (see \cite[\S3]{lazavgit} for further discussion): 
 \begin{itemize}
 \item[(1)] There are finitely many possibilities for the GIT quotients $X\gquot_\calL G$ as one varies the linearization $\calL$. The set of linearizations is partitioned into rational polyhedral chambers parameterizing GIT equivalent linearizations.
 \item[(2)]  The semistable loci satisfy a semi-continuity property. This property induces morphisms between quotients for nearby linearizations.
 \item[(3)] The birational change of the GIT quotient as the linearization moves from one chamber to another by passing a wall is flip like, and can be described quite explicitly. 
 \end{itemize}
The above properties lead to one of main strengths of the theory of Variation of GIT quotients (VGIT): {\it it might be possible to interpolate from an easily understood space $\calM_0$ to a geometrically relevant space $\calM_1$ by varying the linearization in a VGIT set-up}. A spectacular application of this principle is Thaddeus' work \cite{thaddeus0} on Verlinde formula.  A more modest application (but closer to the spirit of these notes) of VGIT  to moduli is Theorem \ref{thmg4git} below. 

\begin{example} We discuss here a simple example of VGIT (for details see \cite[\S3.4]{lazavgit}) which illustrates the interpolation between stability conditions as the linearization varies. Specifically, we consider GIT for pairs $(C,L)$ consisting of a plane cubic and a line.  In this situation, the $G$-linearizations are parameterized by $t\in \bQ$, and they are effective (i.e. there exist invariant sections) for $t\in[0,\frac{3}{2}]$. The GIT stability for the pair $(C,L)$ for parameter $t\in[0,\frac{3}{2}]$ is described as follows:  {\it If  $L$ passes through a singular point of $C$, then the pair $(C,L)$ is $t$-unstable for all $t>0$. Otherwise, $(C,L)$ is $t$-(semi)stable for an interval
 $t\in (\alpha,\beta)$  (resp. $t\in [\alpha,\beta]$), where \begin{equation*}
\alpha=\begin{cases}
0           &\textrm{ if }C \textrm{ has at worst nodes}\\
\frac{3}{5} &\textrm{ if }C \textrm{ has an }A_2 \textrm{ singularity}\\
1           &\textrm{ if }C \textrm{ has an }A_3 \textrm{ singularity}\\
\frac{3}{2} &\textrm{ if }C \textrm{ has a }D_4 \textrm{ singularity}                  
\end{cases}
\textrm{ and }
\beta=\begin{cases}
\frac{3}{5} &\textrm{ if }L \textrm{ is inflectional to }C\\
1           &\textrm{ if }L \textrm{ is tangent to }C\\
\frac{3}{2} &\textrm{ if }L \textrm{ is transversal to }C                  
\end{cases}.
\end{equation*}}
In other words, for $t=0$ the semistability of $(C,L)$ is equivalent to the semistability of $C$; for $t=\frac{3}{2}$ the semistability of $(C,L)$ is equivalent to the semistability of $C\cap L\subset L\cong \bP^1$. The intermediate stability conditions are interpolations of these two extremal conditions.
\end{example}
\subsubsection{Some of the main tools of GIT are: the numerical criterion, Luna's slice theorem, and Kirwan's desingularization.}
In general, it is very hard to describe the (semi)stability conditions for a GIT quotient $X\gquot G$. Essentially, the only effective tool for this is the numerical criterion.  The main points that lead to the numerical criterion are as follows: By definition, the (semi)stability of $x\in X$ is related to the study of orbit closures. By the valuative criterion of properness, to test that an orbit is closed it suffices to consider $1$-parameter families with all fibers in the same orbit. The key point now is that, in the GIT situation, it suffices to study to study $1$-parameter families induced by  $1$-parameter subgroups $\lambda(t)$, $t\in\bC^*$. To such a family, one associates a numerical function $\mu^\mathcal{L}(x,\lambda)$ (where $\calL$ is a $G$-linearization on $X$) which measures the limiting behavior as $t\to0$. This leads to the following:
\begin{theorem}[Hilbert-Mumford Numerical Criterion]\label{numcriterion}
Let $\mathcal{L}$ be an ample $G$-linearized line bundle. Then $x\in X$ is stable 
(resp. semistable) with respect to $\mathcal{L}$ if and only if  $\mu^\mathcal{L}(x,\lambda)>0$ 
(resp. $\mu^\mathcal{L}(x,\lambda)\ge0$) for every nontrivial $1$-PS $\lambda$ of $G$. 
\end{theorem}

The numerical criterion splits the a priori intractable problem of deciding semistability (i.e. finding non-vanishing $G$-invariant sections $\sigma$ at $x\in X$) into two somewhat accessible steps. The first step is a purely combinatorial step involving the weights of the maximal torus in $G$ on a certain representation. Furthermore, it is possible to include the variation of linearization in this combinatorial analysis.  This first step can be effectively solved with the help of a computer (see \cite[\S4.1]{lazavgit}). The second step consists in interpreting geometrically the results of the first step. This is a case by case delicate analysis, and the true bottleneck for the wide applicability of GIT. In any case, for any concrete GIT problem, the numerical criterion gives an algorithmic approach to semistability (see \S\ref{exgit} for a survey of the known GIT examples). In fact, it is reasonable to say that GIT is the most accessible/computable approach to a moduli spaces. 

\smallskip

By construction the quotient $X\gquot_\calL G$ is a normal projective variety. It is natural to ask about its local structure. Of course, the local structure at an orbit $G\cdot x\in X\gquot G$ depends on the local structure of $x\in X$ and the stabilizer $G_x$. For instance if $X$ is smooth, then the geometric quotient $X^s/G$ has only finite quotient singularities (or equivalently, it is a smooth Deligne--Mumford stack) and thus it is well behaved.  In general, for $x\in X^{ss}\setminus X^s$ (wlog assume $G\cdot x$ is a closed orbit, and thus $G_x$ reductive by Matsushima criterion) Luna's theorem gives a precise description of the local structure of the quotient near the orbit $G\cdot x\in X\gquot G$:

\begin{theorem}[Luna Slice Theorem]\label{lunathm}
Given $x\in X^{ss}$ with closed orbit $G\cdot x$ and $X$ smooth,  there exits a $G_x$-invariant normal slice $V_x\subset X^{ss}$ (smooth and affine) to $G\cdot x$ 
such that we have the following commutative diagram with Cartesian squares: 
$$
\begin{CD}
G*_{G_x}\calN_x@<\text{\'etale}<<G*_{G_x}V_x@>\text{\'etale}>>X^{ss}\\
@VVV                               @VVV@VVV\\
\calN_x/ G_x@<\text{\'etale}<<\left(G*_{G_x}V_x\right)/ G@>\text{\'etale}>> X\gquot G\end{CD}
$$
where $\calN_x$ is the fiber at $x$ of the normal bundle to the orbit $G\cdot x$.
\end{theorem}

The main point here is that this theorem reduces (locally) the GIT quotient $X\gquot G$ to a GIT quotient by a smaller reductive group $G_x\subset G$ (typically $G_x$ is a torus, or even $\bC^*$). As an application of this, it is easy to resolve the GIT quotients $X\gquot G$. Specifically, one orders the stabilizers $G_x$ in the obvious way, and then blows-up the strata corresponding to the maximal stabilizers $G_x$. Roughly, (locally) this corresponds to the blow-up of the origin in $\calN_x$. Since the action of $G_x$ is distributed on the exceptional divisor of $\mathrm{Bl}_0\calN_x$, one sees that after such a blow-up the resulting stabilizers are strictly contained in $G_x$. Repeating the process inductively, one arrives to the ideal situation when all the stabilizers are finite.  This resolution process was  developed by Kirwan \cite{kirwan}. As an application, Kirwan used this desingularization process to compute the cohomology of several moduli spaces (e.g. \cite{kirwancoh}). As an aside, we note that Luna's theorem and Kirwan desingularization procedure explain why various moduli spaces are nested. For instance, the exceptional divisor obtained by resolving the ``worst'' singularity for the GIT quotient for cubic fourfolds (\cite{laza1}) is naturally identified with the moduli of degree $2$ $K3$ surfaces (this is also related to the discussion of \S\ref{githodge}). 

\begin{remark}
For a survey of the singularities of GIT quotients see \cite[\S5.1]{lazavgit}. For a discussion of Luna's theorem in connection to VGIT see \cite[\S4.2]{lazavgit}.
\end{remark}
\subsubsection{The main advantage of GIT is that it gives projective models for moduli spaces with weak modular meaning.} Specifically, by construction a GIT quotient $X\gquot G$ is a normal projective variety (assuming, as above, $X$ normal projective). Furthermore, each point of the quotient $X\gquot_\calL G$ corresponds to a unique closed orbit. If there are no strictly semistable points, then $X\gquot G$ is also a geometric quotient, and thus the coarse moduli space of a proper Deligne-Mumford moduli stack. This is typically not the case: there exist strictly semistable points, and thus multiple orbits correspond to the same point in the quotient  $X\gquot G$. It is (typically) not possible to define a functor that selects only the closed orbits. Nonetheless, the weak modular properties of GIT quotients might be the best that one can expect short of a DM stack (see Alper \cite{alpergit} for a formalization of ``good'' moduli stacks). Note also, the following (a consequence of the properness of $X\gquot G$ and of the $1$-to-$1$ correspondence of the points in $X\gquot G$ with the closed orbits):

 \begin{proposition}[GIT Semistable replacement lemma]\label{ssreplace}
Let $S=\Spec R$ and $S^*=\Spec (K)$, where $R$ is a DVR with  field of fractions $K$ and closed point $o$. Assume that $S^*\to X^s/G$ for some GIT quotient. Then, after a finite base change $S'\to S$ (ramified only at the special point $o$), there exists a lift $\tilde{f}:S'\to X^{ss}$ of $f$ as in the diagram:
$$
\xymatrix
{ 
&S'\ar@{->}[r]^{\tilde{f}}\ar@{->}[ld]&X^{ss}\ar@{->}[rd]\\
S&S^*\ar@{_{(}->}[l]\ar@{->}[r]^{f}&X^s/G\ar@{^{(}->}[r]&X\gquot G\\
}.
$$
 Furthermore, one can assume that $\tilde{f}(o)$ belongs to a closed orbit. 
\end{proposition} 

In other words, while a GIT quotient typically fails to have a modular meaning at the boundary, one can use this lemma  to understand the degenerations of smooth objects and then construct or understand a good compactification of the moduli space. A concrete application of this principle is discussed in \S\ref{githodge} below. For some further discussion and some examples see \cite[Ch. 11]{yanosur}.

\subsubsection{Unfortunately, GIT only sees ``linear'' features of the parameterized varieties. Consequently, quite degenerate objects might be semistable, leading to bad singularities for the GIT quotient $X\gquot G$ and the failure of modularity.}\label{disgit} We have mentioned several drawbacks of the GIT approach to moduli spaces: It is difficult to decide stability of objects. Also, typically  there exist strictly semistable points. At these points, the GIT quotient is quite singular, and fails to be modular. The hidden geometric reason for these issues is that GIT only tests for linear features of the objects under consideration (e.g. see Remark \ref{remlogcan}). By considering asymptotic GIT (higher and higher embeddings of a given object) more of the geometric features of the varieties will  be visible by means of ``linear tests''. Unfortunately, as discussed in \S\ref{kstability}, the asymptotic GIT approach is not well behaved/well understood.  

The main point we want to emphasize here is that the unstable objects always satisfy some special conditions with respect to a flag of linear subspaces (e.g. an unstable hypersurface will always contain a singular point and there will be a special tangent direction through this point). This follows for instance from the work of Kempf \cite{kempf}. Namely, for unstable points, there is a distinguished $1$-PS $\lambda$ that destabilizes $x$ (essentially minimizing $\mu(x,\lambda)/|\lambda|$). Then, $\lambda$ determines  a parabolic subgroup $P_\lambda$, which in turn is equivalent to a (partial) flag. The failure of stability involves some special geometric properties with respect to this flag. Consequently, objects that behave well with respect to linear subspaces will tend to be semistable. For instance a conic of multiplicity $\frac{d}{2}$ is semistable when viewed as a curve of degree $d$. Of course, such objects would be disallowed by other more ``modular'' approaches. For instance, for $d=4$ (plane quartics, or genus $3$ curves), the double conic should be replaced by hyperelliptic curves.

\subsection{Applications of GIT to moduli} In this section, we briefly review the scope of GIT constructions in moduli theory. 
\subsubsection{Survey of GIT constructions in moduli theory}\label{exgit} GIT and moduli were tightly connected for over a hundred years. Initially, in late 1800s/early 1900s, the focus was on computing explicitly the rings of invariants for various quotients $X/G$ (for example the ring of invariant polynomials for cubic surfaces). After Hilbert's proof of the finite generation of the ring of invariants $R^G$, the search of explicit invariants fell out of favor. Mumford \cite{GIT} revived GIT  to show that the moduli space of curves $M_g$ is quasi-projective (\cite[Thm. 7.13]{GIT}). The case of abelian varieties (with level structure) is also discussed in Mumford's monograph (\cite[Thm. 7.9]{GIT}).  A little later, Mumford and Gieseker (\cite{mumford}) proved, via GIT, that  the coarse moduli space $\overline{M}_g$ associated to the Deligne--Mumford compactification  of the moduli space of genus $g$ curves is a projective compactification of $M_g$  (for a discussion of this and related constructions, see the survey \cite{msurvey}). Some other major results around the same time include the proof of quasi-projectivity for the moduli of surfaces of general type (Gieseker \cite{gieseker}) and compactifications for the moduli spaces of vector bundles over curves (Mumford, Narasimhan, Seshadri, e.g. \cite{seshadri}) and surfaces (Gieseker \cite{giesekervb}). More recently, Viehweg \cite{viehweg} proved the quasi-projectivity of moduli of varieties of general type (see also \ref{projmodt} and \ref{projmod} below) by using non-standard linearizations on the moduli space. 

The GIT constructions for $\overline{M}_g$  or moduli of surfaces of general type involve {\it asymptotic GIT}, i.e. given a class of polarized varieties $(V,L)$ one considers the GIT quotient of the Chow variety $\mathrm{Chow}_k$ (or Hilbert scheme) for higher and higher embeddings $V\to \bP^{N(k)}$ given by $L^k$ (for $k\gg 0$). In the case of curves the quotients $\mathrm{Chow}_k\gquot \SL(N_k)$ stabilize and give $\overline M_g$ (in fact $k\ge 5$ it is enough).  For surfaces, as discussed in \S\ref{kstability}, there is no stabilization for the asymptotic GIT and it is unclear how to use this asymptotic approach to construct a compact moduli space. 

In recent years, in connection to the Hassett--Keel program, there is a renewed interest in understanding non-asymptotic GIT models for $\overline M_g$. Other non-asymptotic GIT quotients that were studied include moduli for some hypersurfaces:  plane sextics (\cite{shah}), quartic surfaces (\cite{shah4}), cubic threefolds (\cite{allcock}), cubic fourfolds (\cite{laza1}), and some complete intersections (e.g. \cite{avmir}, \cite{g4git}, and \cite{benoist}). We emphasize that for hypersurfaces it is possible to give in an algorithmic way the shape of equations defining unstable hypersurfaces. However, it is difficult to interpret geometrically  the stability conditions. In fact, the higher the degree the worst the singularities that are allowed for stable objects. Consequently, GIT will give a somewhat random compactification for the moduli of  hypersurfaces.

\subsubsection{GIT and the Hassett--Keel program}\label{hkprog} It is of fundamental interest to understand the birational geometry of $\overline M_g$ (see the survey \cite{farkas}), in particular the canonical model $\overline M^{can}_g$ (for $g\ge 24$).  A fundamental insight (due to Hassett and Keel) says that one can approach this problem via interpolation (see esp. \cite{hh,hh2}). Namely, one defines
$$\overline{M}_g(\alpha)=\Proj(R(\overline{M}_g, K_{\overline{M}_g}+\alpha \Delta)),$$
where $\Delta$ is the boundary divisor in $\overline{M}_g$. For $\alpha=1$ (and all $g$) one gets the Deligne-Mumford model $\overline{M}_g$, while for $g\ge 24$ and $\alpha=0$ one gets $\overline M^{can}_g$. It was observed that some of the $\overline{M}_g(\alpha)$ models have modular meaning. For instance, $\overline{M}_g(\frac{9}{11})$ is a moduli of pseudo-stable curves, i.e. curves with nodes and cusps and without elliptic tails. It is conjectured that all $\overline M_g(\alpha)$ for $\alpha\in[0,1]$ have some (weak) modular meaning, and that they behave similarly to the spaces in a VGIT set-up. In other words, (conjecturally) via a finite number of explicit and geometrically meaningful modifications one can pass from $\overline{M}_g=\overline{M}_g(1)$ to $\overline M^{can}_g=\overline M_g(0)$. In this way, one would obtain a satisfactory description of the canonical model $\overline M^{can}_g$ as well as a wealth of information on the birational geometry of $\overline{M}_g$.

Most of the constructed $\overline{M}_g(\alpha)$ spaces were obtained via GIT (see \cite{fs} and \cite{ah} for some recent surveys). Roughly speaking, GIT tends to give compactifications with small boundary. Consequently, once a GIT model $\overline M_g^{GIT}$ with the correct polarization was constructed, one can show that it agrees with $\overline{M}_g(\alpha)$  on open subsets with high codimension boundaries in $\overline{M}_g(\alpha)$  and $\overline M_g^{GIT}$ respectively (and thus $\overline{M}_g(\alpha)$ and $\overline M_g^{GIT}$  agree everywhere). The important point here is that $\overline M_g^{GIT}$ is a projective variety and that its polarization can be easily understood by descent from a parameter space (typically a Chow variety or Hilbert scheme for small embeddings of genus $g$ curves). Another important point is that  $\overline M_g^{GIT}$ comes with a (weak) modular interpretation and thus it induces such a modular interpretation for $\overline{M}_g(\alpha)$. 

As already mentioned,  the behavior of $\overline M_g(\alpha)$ seems to be parallel to that of quotients in a VGIT set-up. One might conjecture that there is a master VGIT problem modeling all of $\overline M_g(\alpha)$. Unfortunately, the best result so far is: 

\begin{theorem}[\cite{g4git}]\label{thmg4git}
For $\alpha \le \frac{5}{9}$, the log minimal models $\overline{M}_4(\alpha)$ arise as GIT quotients of the parameter space $\bP E$ for $(2,3)$ complete intersections in $\bP^3$.  Moreover, the VGIT problem gives us the following diagram:

\begin{equation}\label{mainthm}
\xymatrix @R=.07in @C=.07in{
 & \overline{M}_4 ( \frac{5}{9} , \frac{23}{44} ) \ar@{-->}[rr]\ar[ldd] \ar[rdd] & & \overline{M}_4 ( \frac{23}{44} , \frac{1}{2} ) \ar[ldd] \ar[rdd] \ar@{-->}[rr]& & \overline{M}_4 ( \frac{1}{2} , \frac{29}{60} ) \ar[ldd] \ar@{=>}[rdd] & \\
&&&&&& \\
\overline{M}_4 ( \frac{5}{9} ) & & \overline{M}_4 ( \frac{23}{44} ) & & \overline{M}_4 ( \frac{1}{2} ) & & \overline{M}_4 [\frac{29}{60},\frac{8}{17}) \ar[dd] \\
&&&&&&\\
 & & & & & & \overline{M}_4 (\frac{8}{17} )   = \{*\}
 }
 \end{equation}

 More specifically,
 \begin{itemize}
 \item[i)] the end point $\overline{M}_4 ( \frac{8}{17}+\epsilon)$ is obtained via GIT for $(3,3)$ curves on $\bP^1 \times \bP^1$ as discussed in \cite{maksymg4};
 \item[ii)] the other end point $\overline{M}_4 ( \frac{5}{9} )$ is obtained via GIT for the Chow variety of genus 4 canonical curves as discussed in \cite{g4ball};
 \item[iii)] the remaining spaces   $\overline{M}_4 ( \alpha ) $ for $\alpha$ in the range $\frac{8}{17} < \alpha<  \frac{5}{9}$  are obtained via appropriate $\Hilb^m_{4,1}$ quotients, with the exception of $\alpha= \frac{23}{44}$.
 \end{itemize}
\end{theorem}

\begin{remark}
The  Hassett--Keel program is currently established for either values of $\alpha$ close to $1$ (e.g. \cite{hh,hh2}) or small values of the genus $g$ (e.g. \cite{hl}, \cite{maksymg4}, \cite{g4ball,g4git}, \cite{fs5}).  For some general predictions on $M_g(\alpha)$ see \cite{alper}. 
\end{remark}

\section{Moduli and periods}\label{sectht}
A different approach to the construction of moduli spaces is based on the idea of associating to a variety its cohomology, and studying the induced variation of Hodge structures (VHS). As we will discuss here, the scope of this approach is quite limited in practice. However, when the Hodge theoretic approach is applicable, it has strong implications on the structure of the moduli space; and thus this is a highly desirable situation. We discuss some examples of moduli space for which both the GIT and Hodge theoretic approach are applicable. Each approach gives a different facet of the moduli space. 

Some general references for the material discussed here include \cite{carlson}, \cite{voisin}, \cite{ggk}, \cite{milne}, and \cite[Ch. 2, 3]{ueno}.

\subsection{Period maps}
The primitive cohomology of a smooth projective variety carries a polarized Hodge structure such that the associated Hodge filtration varies holomorphically in families. 
More formally, one says that for a smooth family $\pi:\calX\to S$ of algebraic varieties, $(R^n\pi_*\underline{\bZ}_\calX)_{prim.}$ defines a polarized {Variation of Hodge Structures} (VHS) over $S$, which in turn defines a period map $\calP:S\to \bD/\Gamma$. In order to use $\calP$ to construct a moduli space, we need to discuss the injectivity of $\calP$ (``Torelli Theorems'') and the image of $\calP$.

\subsubsection{Period domains and period maps.} The period domain $\bD$ is the classifying space of Hodge structures of a given type. Specifically, the polarized Hodge structures of weight $n$ satisfy the Hodge-Riemman bilinear relations:
\begin{itemize}
\item[(HR1)] $F^p=(F^{n-p+1})^\perp$;
\item[(HR2)] $(-1)^{n(n-1)} i^{p-q}(\alpha,\bar\alpha)>0$ for $\alpha \in H^{p,q}=F^p\cap \overline{F^q}$ (with $p+q=n$).
\end{itemize}
 The first condition (HR1) defines a projective homogeneous variety $\check \bD=\bG_\bC/\bB$, a subvariety of  a flag manifold. The condition (HR2) gives that $\bD$  is an open  subset (in the classical topology) of $\check \bD$, in particular a complex manifold. Note also that the period domain is homogeneous $\bD=\bG_\bR/\bK$ (with $\bK=\bB\cap \bG_\bR$ a compact subgroup of $\bG_\bR$) and {\it semi-algebraic} (given by algebraic inequalities involving $\mathrm{Re}$ and $\mathrm{Im}$ of holomorphic coordinates). It is important to note that there are only two cases when $\bD$ is a  Hermitian symmetric domain (or equivalently $\bK$ is a maximal compact subgroup):  
\begin{itemize}
\item weight $1$ Hodge structures (abelian variety type): $\bD$ is the Siegel upper half space $\mathfrak{H}_g=\{ A\in M_{g\times g}(\bC)\mid A=A^t, \mathrm{Im}(A)>0\}\cong \Sp(2g)/U(g)$;
\item weight $2$ Hodge structures with $h^{2,0}=1$ ($K3$ type):  $\bD$ is a Type IV domain $\{\omega\in \bP(\Lambda_\bC)\mid \omega.\omega=0, \omega.\bar\omega>0\}\cong \SO(2,n)/S(O(2)\times O(n))$, where $\Lambda$ is a lattice of signature $(2,k)$ (here $\check\bD$ is a quadric hypersurface in $\bP^{k+1}$; thus $\dim \bD=\dim \check\bD=k$).
\end{itemize}

\begin{example}
For Hodge structures of Calabi-Yau threefold type with Hodge numbers $(1,h,h,1)$, the period domain is 
$$\bD=\Sp(2(1+h))/U(1)\times U(h).$$
The maximal compact subgroup is $U(1+h)$ and the inclusion $U(1)\times U(h)\subset U(1+h)$ induces a natural map $\bD\to \mathfrak{H}_{1+h}$ which is neither holomorphic or anti-holomorphic (only real analytic). 
\end{example}

A {\it Variation of Hodge Structure} over $S$  is a triple $(\calV,F^\bullet,\nabla)$ consisting of a flat vector bundle $(\calV,\nabla)$ together with holomorphic subbundles $F^n\subset F^{n-1}\dots\subset F^0=\calV$. We assume that the VHS is polarized (i.e. there is a compatible bilinear form such that (HR1) and (HR2) are satisfied). By passing to a trivialization of the associated local system we obtain a period map:
$$\widetilde P:\widetilde S\to \bD,$$
where $\widetilde S$ is the universal cover of $S$.  Alternatively, we get
$$\calP:S\to \bD/\Gamma$$
where $\Gamma=\mathrm{Im}(\pi_1(S)\to \bG_Z)$. Note that $\calP$ is a locally liftable analytic map (N.B. $\bG_Z$, and thus $\Gamma$, acts properly discontinuous on $\bD$, making $\bD/\Gamma$ an analytic space). 

\subsubsection{Griffiths' Transversality. Classical and semi-classical period maps.}The period maps arising from algebraic geometry satisfy an additional condition: Griffiths transversality 
$$\nabla F^{p}\subseteq F^{p-1}\otimes \Omega_S$$
which is a non-trivial condition except for the case when $\bD$ is Hermitian symmetric. As noted above, there are only two cases when $\bD$ is Hermitian symmetric: abelian varieties and $K3$ type. 

We are interested in the situation when $\calM$ is a moduli space (stack) of smooth varieties, and 
$$\calP:\calM\to \bD/\Gamma$$
is induced by the variation of cohomology $R^n\pi_*\underline \bZ$. For constructing a moduli space via periods, we would like that $\calP$ is a birational map (a priori the period maps are analytic, but there are various algebraicity results; we ignore the issue here).  There are two statements that one tries to prove: $\calP$ is injective and $\calP$ is dominant. Due to the Griffiths transversality, except for ppavs and K3 type,  $\calP$  is never dominant. Even worse, the image of the period map is typically highly transcendental. Namely, say that $Z$ (a closed analytic subvariety) is the image of a period map in $\bD$.  Then $Z\subset \bD$ is horizontal w.r.t. the distribution given by Griffiths' transversality. Furthermore, in the algebro-geometric situation,  $Z$ is stabilized by a big subgroup $\Gamma\subset G_\bZ$. In fact, it is reasonable to assume $Z/\Gamma$ is quasi-projective. Then, we proved the following:
\begin{theorem}[{\cite{cypaper}}]\label{thmalg} Let $Z$ be a closed horizontal subvariety of a classifying space $\bD=\bG_\bR/\mathbf K$ for Hodge structures and let $\Gamma$ be the stabilizer of $Z$ in the appropriate arithmetic group $\bG_\bZ$. Assume that 
\begin{itemize}
\item[(i)] $S=\Gamma\backslash Z$ is strongly quasi-projective;
\item[(ii)] $Z$ is semi-algebraic in $\bD$ (i.e. open in its Zariski closure in $\check \bD$). 
\end{itemize}
Then  $Z$ is a Hermitian symmetric domain  whose embedding  in $\bD$ is an equivariant, holomorphic, horizontal embedding.
\end{theorem} 

In other words, the only cases when the image of a period map can be described purely algebraically are the Shimura type cases (see \cite{milne} and \cite{dshimura}). They are slight generalizations of the classical cases of  ppavs and K3 type;  we call them {\it semi-classical cases}. A number of such semi-classical examples are discussed in the recent literature (see for instance 
 the series of ball quotient examples of Kondo \cite{dk}). These semi-classical examples are naturally understood in the context of Mumford-Tate (MT) domains \cite{ggk}. Essentially, a MT domain is the smallest homogeneous subdomains of $\bD$ defined over $\bQ$ which contains the image of a VHS. To belong to a MT subdomain of $\bD$ is equivalent to saying that there exist some special Hodge tensors for the VHS. The simplest situation is that of Hodge structures with extra endomorphisms as in the following example:
\begin{example}[Kondo \cite{kondo3}]
A generic genus $3$ curve $C$ is a plane quartic. To it one can associate a quartic $K3$ surface $S$ by taking the $\mu_4$-cover of $\bP^2$ branched along $C$. For example, if $C=V(f_4)$, then $S=V(f_3(x_0,x_1.x_2)+x_3^4)\subset \bP^3$. In this way, one gets a period map
$$\calP:M_3^{nh}\to \calF_4\cong\calD/\Gamma$$ 
from the moduli of non-hyperelliptic genus $3$ curves to the period domain of degree $4$ $K3$ surfaces. Since the resulting $K3$ Hodge structures are special (they have multiplication by $\mu_4$), the image of $\calP$ will land into a MT subdomain. In this situation, the MT subdomain will be a $9$-dimensional complex ball $\calB$ embedded geodesically into the $19$-dimensional Type IV domain $\calD$. In conclusion, one gets 
$$\calP:M_3^{nh}\to \calB/\Gamma'\subset\calD/\Gamma,$$
which turns out to be birational (see \cite{l3} and \S\ref{githodge}). This construction is also relevant to the Hassett--Keel program in genus $3$ (e.g. \cite{hl}).  
\end{example}

We emphasize that the situation covered by Theorem \ref{thmalg} is very special. For instance, while there exist moduli of Calabi-Yau threefolds that have period maps to Hermitian symmetric domains as in the theorem (e.g. the Borcea--Voisin \cite{borcea} and Rohde--van Geemen \cite{rohde} examples), most moduli of Calabi-Yau threefolds (e.g. for quintic threefolds) are not of this type (see \cite{cypaper}).

\subsubsection{Torelli theorems} The differential of a period map $\calP:S\to \bD/\Gamma$ for a family of algebraic varieties has a simple description in cohomological terms.  Namely, since the tangent space to a Grassmannian at the point $F\subset H_\bC$ is canonically $\Hom(F,H_\bC/F)$,  it follows that the tangent space to a period domain satisfies $T_o\bD\subseteq \oplus_p\Hom(F^p,H_\bC/F^p)$. Due to the Griffiths' transversality, we then view the differential of  $\calP$ as a map
$$d\calP:T_{S,s}\to \oplus_p \Hom(F^p/F^{p+1},F^{p-1}/F^{p})\subset T_{\calP(s)}\bD.$$
 Now $F^{p}/F^{p+1}=H^{p,q}=H^q(X_s,\Omega_{X_s}^p)$. Then,  the differential of the period map is given by (see \cite[I \S10.2.3]{voisin}):
 \begin{equation*}d\calP(\xi)(\phi)=\kappa(\psi)\cup\phi \in H^{q+1}(X_s,\Omega_{X_s}^{p-1}) \textrm{ for } \xi\in T_{S,s} \textrm{ and } \phi\in H^q(X_s,\Omega_{X_s}^p),\end{equation*}
 where $\kappa:T_{S,s} \to H^1(X_s,T_s)$ is the Kodaira-Spencer map, and $$\cup:H^1(X_s,T_s)\times H^q(X_s,\Omega_{X_s}^p)\to H^{q+1}(X_s,\Omega_{X_s}^{p-1})$$ is induced by the contraction map $T_s\otimes \Omega_{X_s}^p\to \Omega_{X_s}^{p-1}$. As a consequence, it is not hard to prove {\it infinitesimal Torelli} in a variety of cases.  An interesting case is that of Calabi-Yau threefolds. Due to the Bogomolov-Tian Theorem, the moduli space $\calM$ of Calabi-Yau threefolds (in a fixed deformation class)  is smooth  and 
 $$T_s\calM\cong H^1(X_s,T_s)\cong H^1(X_s,\Omega^2_{X_s})$$
 (using the Calabi-Yau condition: $K_{X_s}$ is trivial). Since $H^{3,0}(X_s)\cong \bC$, it follows that $d\calP$ gives an isomorphism between the tangent space to moduli and the subspace of horizontal directions in $T_{\calP(s)}\bD$. Thus, locally $\calM$ is identified with a maximal horizontal subvariety $Z\subset \bD$ (N.B. $\dim \calM=\dim Z= h^{2,1}$). A local description of the horizontal subvarieties $Z$ of $\bD$ in the Calabi-Yau case was given by Bryant--Griffiths \cite{bg}. A global description is much harder, and little is known (see however \cite{morgan} and \cite{cypaper}).
 
Global Torelli results ($\calP$ injective) are harder to obtain (as they involve the global geometry of the moduli spaces). In particular, the global Torelli is known to hold in the following cases: 
\begin{itemize}
\item for abelian varieties (this is essentially a tautological statement);
\item for curves  (a classical statement), and cubic threefolds (\cite{cg});
\item for $K3$s  (a nontrivial result, typically proved using the density of the Kummer surfaces in the moduli of $K3$s, and then applying Torelli for ppavs);
\item for $K3$ like situations: e.g.  cubic fourfolds (\cite{voisincubic}) and compact Hyperk\"ahler manifolds (due to Verbitsky, see \cite{huybrechtsv}). 
\end{itemize}
It is also known that generic Torelli ($\calP$ is generically injective) theorem holds for most hypersurfaces (a result due to Donagi, see \cite[II \S6.3.2]{voisin}). Roughly, the idea here is that the image of the period map has such a  high codimension (satisfies so many infinitesimal conditions), that the hypersurfaces can be recovered from infinitesimal data. Of course, the global or even generic Torelli doesn't always hold. For example, for del Pezzos the associated VHS is trivial; there are also some non-trivial counterexamples to Torelli. In any case, it is expected that various forms of Torelli will hold quite generally. On the other hand, with few exceptions, the period map is not dominant. 

\subsubsection{For moduli constructions, we restrict to the case $\bD$ is a Hermitian symmetric domain (or more generally to the cases covered by Theorem \ref{thmalg}).} In conclusion, our understanding of the images of period maps outside the situation covered by Theorem \ref{thmalg} is quite limited. Consequently, the only case that we know how to use the period map for constructing a moduli space is when we have a period map 
$$\calP:\calM\to \calD/\Gamma,$$
where $\calD$ is a Hermitian symmetric domain and $\Gamma$ an arithmetic group acting on $\calD$. Due to  Baily-Borel and Borel theorems discussed below, $\calP$ is in fact an algebraic (typically birational) map between quasi-projective varieties. While rare, this situation is particularly good 
due to the special structure of $\calD/\Gamma$; some applications to moduli are explained below.

\begin{remark}
We note here that there is at least one application to using the period map  in the case $\calD$ is not Hermitian symmetric. Namely, the period domains $\bD$ are negatively curved in the horizontal directions. Consequently, if the Torelli theorem holds, one obtains hyperbolicity results for the moduli space (see \cite{kebekus} for a survey on hyperbolicity). Even when Torelli doesn't hold, one might obtain hyperbolicity results by using cyclic cover constructions and applying Torelli theorems on related moduli spaces (e.g. see \cite{vz}). 
\end{remark}

\subsection{Applications of locally symmetric varieties} 
The main advantage of realizing a moduli space as a locally symmetric variety $\calD/\Gamma$ (with $\calD$ a Hermitian symmetric domain and $\Gamma$ an arithmetic group acting on $\calD$) is that such a variety has a lot of additional structure, which brings in additional tools such as the theory of automorphic forms. We briefly survey two applications of the description of moduli as locally symmetric varieties: the existence of natural compactifications and results on the Kodaira dimension of the moduli spaces. 
\subsubsection{Compactifications of locally symmetric varieties}\label{compactdg} 
A locally symmetric variety has several natural compactifications. First, any locally symmetric variety $\calD/\Gamma$ has a canonical compactification, the Satake-Bailly-Borel compactification, which can be defined by:
\begin{equation}\label{defBB}
\calD/\Gamma\subseteq (\calD/\Gamma)^*:=\Proj A(\Gamma),
\end{equation}
where $A(\Gamma)$ is the ring of $\Gamma$-automorphic forms on $\calD$ (this is a finitely generated ring by \cite{bb}).  By construction $(\calD/\Gamma)^*$ is a projective variety which does not depend on any choices. It is also a minimal compactification, in the sense that any normal crossing compactification of $\calD/\Gamma$ will map to it. More precisely, the following holds (due to curvature properties of period domains):
\begin{theorem}[Borel's Extension Theorem \cite{borel}]\label{borelext}
Let $S$ be a smooth variety, and $\overline S$ a smooth simple normal crossing (partial) compactification of $S$. Then any locally liftable map $S\to \calD/\Gamma$ extends to a regular map $\overline S\to (\calD/\Gamma)^*$.
\end{theorem}

The main disadvantage of the SBB  compactification $(\calD/\Gamma)^*$ is that it is quite small, and thus it doesn't accurately reflects the geometry of degenerations. Also, $(\calD/\Gamma)^*$ tends to be quite singular. For instance, it is well known that  
$$\calA_g^*=\calA_g\sqcup \calA_{g-1} \dots \sqcup \calA_0$$
and thus the boundary has codimension $g$. Similarly, for Type IV domains (i.e.  $K3$ type period domains) the boundary is $1$-dimensional. 

To rectify this issue, Mumford et al. \cite{mumfordtai} have introduced the {\it toroidal compactifications} $\overline{\calD/\Gamma}^\Sigma$. By construction $\overline{\calD/\Gamma}^\Sigma$ depends on a choice, an admissible rational polyhedral decomposition $\Sigma$ of a certain cone (see \cite[Def. 7.3]{namikawabook}). The toroidal compactifications come equipped with  natural forgetful maps 
$$\overline{\calD/\Gamma}^\Sigma\to (\calD/\Gamma)^*$$
 for any $\Sigma$. For suitable $\Sigma$ the compactifications $\overline{\calD/\Gamma}^\Sigma$  are smooth (up to finite quotients) and projective. Of course, the main disadvantage of the toroidal compactifications $\overline{\calD/\Gamma}^\Sigma$ is that there is a plethora of choices. It is unclear which of the toroidal compactifications should have a geometric meaning. The following are some known facts on the modular meaning of some toroidal compactifications:
\begin{itemize}
\item[(A)] (Mumford, Namikawa \cite{namikawabook}, Alexeev \cite{alexeevab}) {\it The second Voronoi compactification $\overline\calA^{Vor}_g$ for $\calA_g\cong\mathfrak{H}_g/\Sp(2g,\bZ)$ has a modular interpretation} (in a strong functorial sense).
\item[(C)] (Mumford, Alexeev--Brunyate \cite{ab}) {\it The period map $M_g\to \calA_g$ extends to a regular $\overline{M}_g\to \overline{\calA}_g^{\Sigma}$, where $\Sigma$ is the second Voronoi or perfect cone. It does not extend to the central cone compactification (for $g\ge 9$)}. 
\item[(P)] (Friedman--Smith \cite{fsp}, see also \cite{abh} and \cite{cmghl}) {\it The period map $R_g\to \calA_g$, where $R_g$ is the moduli of Prym curves, does not extend to a regular map $\overline{R}_g\to \overline{\calA}_g^{\Sigma}$ for any of the standard toroidal compactifications} (here $\overline{R}_g$ is Beauville's admissible covers compactification of $R_g$). 
\end{itemize}
For moduli of $K3$s, it is an open question to give geometric meaning to any of the toroidal compactifications (analogue to (A) above); the best result so far is Theorem \ref{thmlooijenga} below. Another interesting situation is that of cubic threefolds (which are closely related to genus $4$ curves and Prym varieties of genus $5$). Specifically, the intermediate Jacobian for a cubic threefold is a ppav of genus $5$. One is interested in understanding the closure of the intermediate Jacobian locus in some toroidal compactification $\overline\calA_5^\Sigma$ (this question is related to (C) and (P) above); see \cite{cml} and \cite{gh} for some results on this topic. In particular, in \cite{cml} we have computed the closure of the intermediate Jacobian locus in the SBB compactification $\calA_5^*$. The key ingredients for this result are the extension Theorem \ref{borelext} and an appropriate blow-up of the GIT compactification of the moduli space of cubic threefolds. As discussed in the following subsection, it is quite typical to use auxiliary compactifications (such as GIT) in order to understand the geometric meaning at the boundary of the period domain. 

\begin{remark}
The limit of a degeneration of Hodge structures is a {\it mixed Hodge structure} $(H,F^\bullet, W_\bullet)$. From this point of view, it is known that the data encoded in the boundary of the SBB compactification $(\calD/\Gamma)^*$ is equivalent to the data of the graded pieces $\Gr^W_\bullet H$ of the mixed Hodge structure. For instance, in a degeneration of abelian varieties the limit is a semiabelian variety (i.e. an extension $0\to (\bC^*)^n\to X\to A\to 0$ of an abelian variety $A$ by a torus). The SBB compactification remembers only the compact part $A\in \calA_{g-n}\subset \partial \calA_g^*$. On the other hand, the toroidal compactification $\overline{\calD/\Gamma}^\Sigma$  encodes the full limit mixed Hodge structure (see \cite{cattani}). This gives a conceptual explanation of why (A) is plausible. 
\end{remark}

\begin{remark}
There exist extension theorems from normal crossing compactifications $S\to \overline{S}$ to toroidal compactifications (e.g. \cite[Thm.~7.29]{namikawabook}), but in contrast to   Theorem \ref{borelext},  the extension $\overline S\to\overline{\calD/\Gamma}^\Sigma$ exists only if the cones spanned by the monodromies around the boundary divisors in $\overline S$ are compatible with the cones of the decomposition $\Sigma$ (see also \cite{cattani}). This is due to the fact that $\overline{\calD/\Gamma}^\Sigma$ has a toric structure near the boundary. 
\end{remark}
 
 \begin{remark}
 There are two directions in which the SBB and toroidal compactifications can be generalized. First, for the Hermitian symmetric case, Looijenga \cite{lc1,lc2} has introduced the {\it semi-toric compactifications}. Roughly, one makes a choice $\Sigma$ similar to the toroidal case, giving a compactification $\overline{\calD/\Gamma}^\Sigma$, but (in contrast to the toroidal case) $\Sigma$ is only a {\it locally} rational polyhedral decomposition. The semi-toric construction generalizes both the Baily-Borel (roughly $\Sigma$ is the trivial decomposition) and the toroidal constructions ($\Sigma$ is rational polyhedral).   For an application to moduli of the semi-toric compactifications see \S\ref{githodge}. In a different direction, one can ask for compactifications when $\calD$ is not a Hermitian symmetric domain, but rather a general period domain $\bD$. This case was analyzed by Kato-Nakayama-Usui \cite{katousui}. In this situation one doesn't get a compact space or even a variety, but  a ``log variety''; still this is quite useful in applications (e.g. regarding the algebraicity of various Hodge theoretic loci). 
 \end{remark}

\subsubsection{Using automorphic forms to obtain geometric consequences}
As discussed in the previous subsection, the SBB compactification $(\calD/\Gamma)^*$ is a projective variety, being the $\Proj$ of the finitely generated algebra of automorphic forms. In fact, the automorphic forms can be regarded as sections of a certain ample line bundle. If the arithmetic group $\Gamma$ acts freely (which never happens in practice), the automorphic forms would be pluricanonical forms on $\calD/\Gamma$. Thus, it is no surprise that one can use automorphic forms to prove that various moduli spaces are of general type. In practice, since $\Gamma$ doesn't act freely on $\calD$, one has to take care of the ramification of the natural projection $\calD\to\calD/\Gamma$ and of the singularities at the boundary of  $(\calD/\Gamma)^*$  (typically one has to pass to a toroidal compactification). The situation for ppavs is well understood. In particular, 

\begin{theorem}[Tai \cite{tai}, Mumford \cite{mumag}]
$\calA_g$ is of general type for $g>6$. 
\end{theorem}
For small genus ($g\le 5$), since $\calA_g$ is closely related to $M_g$ or $R_g$ and these have simple GIT models, one obtains unirationality (the case $g=6$ is open).  Similar results hold also for $K3$s. Here, it is harder to produce automorphic forms. Essentially, all automorphic forms that occur in the study of the moduli spaces of $K3$s are obtained by appropriate restrictions of the Borcherds' automorphic form $\Phi_{12}$ (see \cite{bkpsb}). In particular, Gritsenko--Hulek--Sankaran \cite{ghs} have proved that for $d$ large enough, $\calF_d$ is of general type (see also \cite{ghs2} for a survey of related results). 

The are numerous other geometric applications involving automorphic forms. For instance, we have used a restriction of the $\Phi_{12}$ form to relate Kondo's ball quotient model for $M_4$ to a step in the Hassett-Keel program (this is related to Theorem \ref{thmg4git} and \S\ref{githodge} below).  Maulik--Pandharipande \cite{mp} have given a generating function for the degrees of the Noether-Lefchetz  loci in the moduli of K3 surfaces with applications to enumerative geometry (see also \cite{lz} for the similar case of cubic fourfolds). Even more spectacularly, Maulik and Charles \cite{charles} have used related ideas to give a proof of the Tate conjecture for $K3$s. 

In conclusion, we see that a (birational) model of type $\calD/\Gamma$ for a moduli space has deep consequences; essentially the global results obtained by this construction are unparalleled. As already explained, the main drawback is that there is limited applicability of this construction. Another issue is that it is typical difficult to understand the geometric meaning of the compactifications $(\calD/\Gamma)^*$ and $\overline{\calD/\Gamma}^\Sigma$. To handle this last issue, one typically needs to construct a different (more geometric) compactification $\overline\calM$ for the moduli space, and then to study the behavior of the period map at the boundary of $\overline\calM$. We discuss a situation  when this comparison approach works in the following subsection (see also \S\ref{sdubois}).  

\subsection{Comparison to GIT compactifications}\label{githodge}
As explained above, the ideal situation is to have a period map 
$$\calP^\circ:\calM^\circ\to \calD/\Gamma,$$
from a moduli space $\calM^\circ$ of smooth varieties to the quotient of a Hermitian symmetric domain $\calD$ by an arithmetic group $\Gamma$. Even after assuming that the Torelli theorem holds and that $\calP$ is dominant (and thus $\calP$ is a birational morphism of quasi-projective varieties) it is typically hard to understand the image of the period map and the geometric meaning of the boundary (the complement of $\mathrm{Im}(\calP^\circ)$). We describe here a situation when we can answer these questions by using an alternative GIT compactification. As a byproduct, one obtains a dual GIT/Hodge theoretic description for a moduli space. This has the advantage of both being (weakly) geometric due to GIT, and having a lot of additional structure due to the $\calD/\Gamma$ description.

Namely, with assumptions as above, we further restrict to the situation that $\calD$ is either a Type IV domain (i.e. $K3$ type) or a complex ball. What is special about these cases is that the maximal {\it Noether-Lefschetz loci} (i.e. loci where the Hodge structure is special in the sense of Mumford-Tate groups, see \cite[\S IIC]{ggk}) are {\it Heegner divisors}, i.e. codimension $1$ subdomains $\calH\subset \calD$ of the same type as $\calD$ (e.g. for a $n$-dimensional complex ball $\calD$, the loci $\calH$ will be $(n-1)$-dimensional sub-balls in $\calD$). In fact, these Heegner divisors are nothing but hyperplane sections of $\calD$ in an appropriate sense (i.e. using $\calD\subset \check\calD\subseteq \bP^N$). In this situation ($\calP$ birational and $\calD$ of Type IV or ball), Looijenga \cite{lc2,ls} has observed the following:
\begin{itemize}
\item[(i)] $\calP^\circ$ is typically an open embedding with image the 
complement of a $\Gamma$-invariant arrangement of hyperplanes. 
\item[(ii)] Given (i),  $\calP^\circ$ can be (typically) extended to an explicit isomorphism between a GIT compactification of $\calM^\circ$ and a semi-toric modification of $\calD/\Gamma$. 
\end{itemize}

The idea for proving these statements is to consider a GIT compactification $\overline\calM$ of $\calM^\circ$ which almost agrees with $(\calD/\Gamma)^*$ as a polarized variety. More specifically, typically the period map $\calM^\circ\to \calD/\Gamma$ extends naturally to a larger moduli stack $\calM$ of very mildly singular varieties. For example, for polarized $K3$ surfaces the period map extends naturally to the moduli $\calM$ of polarized $K3$s with ADE singularities. Then one constructs a GIT compactification $\overline \calM$ such that $\calM\subset \overline \calM$ has boundary of codimension $2$ and with the right polarization (explained below). The GIT quotients typically satisfy the first condition, as their boundary is quite small (however this is typically not true for the smooth locus $\calM^\circ\subset \overline\calM$). 
Then one considers the birational map 
$$\calP:\overline\calM \dashrightarrow (\calD/\Gamma)^*,$$
which is regular on $\calM$. Using the GIT semistable replacement lemma \ref{ssreplace}, one proves (i) or even 
\begin{equation}\label{isoopen}
\calM\cong (\calD\setminus \calH)/\Gamma
\end{equation}
as quasi-projective varieties. Roughly, one shows that the limit Hodge structure for a geometric degeneration with central fiber in a list (determined by GIT) will be either not pure (and thus will go to the boundary of $\calD/\Gamma$) or it will be Hodge special, and thus it will belong to a Noether-Lefschetz locus (a hyperplane in  $\calD$). 

Having the right polarization for the GIT quotient $\overline \calM$ means that the isomorphism \eqref{isoopen} preserves the line bundles induced from the projective  varieties $\overline \calM$ and $(\calD/\Gamma)^*$ respectively. By construction, we have that $\calM\subset \overline \calM$ has boundary of codimension at least $2$. If the same were true for $(\calD\setminus \calH)/\Gamma\subset (\calD/\Gamma)^*$ (e.g. if $\calH=\emptyset$), then due to the extension of sections of line bundles in codimension $2$, we would get $\overline \calM\cong (\calD/\Gamma)^*$. In general the hyperplane arrangement $\calH$ is nontrivial and thus $\calH/\Gamma$ is a divisor.  Looijenga noted that by using an appropriate semi-toric compactification $\overline{\calD/\Gamma}^\Sigma$ (where the decomposition $\Sigma$ is induced by the arrangement $\calH$), one obtains that the divisor $\calH/\Gamma$ is $\mathbb Q$-Cartier and negative (and thus can be contracted).   By this process, Looijenga obtains a modification $\overline{(\calD/\Gamma)}_{\calH}$ of  $(\calD/\Gamma)^*$ with the property that it preserves $(\calD\setminus \calH)/\Gamma$, but contracts $\calH/\Gamma$ to a smaller dimensional stratum (and thus  the codimension $2$ argument holds).  We conclude
\begin{equation}\label{comparison}
\overline\calM\cong \overline{(\calD/\Gamma)}_{\calH},
\end{equation}
(see \cite[Thm. 7.6]{lc2} for a precise statement and precise assumptions).

\begin{remark}
Similarly to the definition \ref{defBB} for the SBB compactification $(\calD/\Gamma)^*$,  one can describe $\overline{(\calD/\Gamma)}_{\calH}$ as the $\Proj$ of a ring of meromorphic automorphic forms with poles along $\calH$. More explicitly, $\overline{(\calD/\Gamma)}_{\calH}$ is obtained by (1) doing a small partial resolution of the boundary $(\calD/\Gamma)^*$ so that $\calH/\Gamma$ becomes $\bQ$-Cartier, followed by (2)  blowing up the arrangement $\calH$ to normal crossing and finally (3) contracting in the opposite direction (i.e. the arrangement is ``flipped'': linear strata in the arrangement $\calH$ are replaced in a way that changes the dimension with codimension). The main points are that $\overline{(\calD/\Gamma)}_{\calH}$ is easy to describe in practice (e.g. see Theorem \ref{thmlooijenga} below) and that this space has a structure (e.g. various arithmetic stratifications) similar to $(\calD/\Gamma)^*$. 
\end{remark}

The simplest example of \eqref{comparison} is the moduli space of elliptic curves. On one hand, we can give a compactification $\overline \calM$ for this moduli space by using the GIT quotient for plane cubics. On the other hand, the moduli space of elliptic curves can be described as the quotient $\calh/\SL(2,\bZ)$ of the Siegel upper half space by the modular group.  As it is well known, we have 
$$
\overline{ M}_1\cong \bP^1\cong \left(\mathfrak H/\SL(2,\bZ)\right)^*\cong \bP\Sym^3 V^*\gquot\SL(3,\bC).
$$
 Some analogue results, when a GIT quotient $\overline \calM$ is isomorphic to a Baily-Borel compactification $(\calD/\Gamma)^*$ (i.e. the arrangement $\calH=\emptyset$ in \eqref{comparison}) hold for some of the $M_{0,n}$ with $n\le 12$ (Deligne--Mostow \cite{dmo}) and the moduli of cubic surfaces (cf. \cite{act1}, \cite{dgk}). In general, the hyperplane arrangement $\calH$ is non-trivial, and thus the GIT quotient $\overline \calM$ is isomorphic to a  semi-toric modification of $\calD/\Gamma$ (which is a non-trivial explicit modification of $(\calD/\Gamma)^*$). Specifically, the following moduli spaces are known to have a dual GIT/Hodge theoretic description:
\begin{itemize}
\item $M_g$  for $g\le 4$ and $g=6$ (\cite{kondo3}, \cite{l3}, \cite{kondo4}, \cite{g4ball}, \cite{ak});
\item moduli of del Pezzo surfaces (\cite{act1}, \cite{dgk}, \cite{heclo});
\item moduli of low degree $K3$ surfaces (\cite{shah,shah4}, \cite[\S8.2]{lc2});
\item moduli of cubic threefolds (\cite{allcock, act2}, \cite{ls});
\item moduli of cubic fourfolds (\cite{laza1,laza2}, \cite{looijenga}). 
\end{itemize}
The example of degree $2$ $K3$ surfaces is discussed in some detail below. 

\smallskip

We note that the comparison result \eqref{comparison} is quite surprising given the different nature of the objects under consideration: on one hand the GIT quotient is purely algebraic, while $\calD/\Gamma$ is of analytic and arithmetic nature. This explains to  a certain extent the relative scarcity of such examples.

\subsubsection{An explicit example - the moduli of degree $2$ $K3$ surfaces} Let $\calF_2$ be the moduli space of degree $2$ $K3$ surfaces. As it is well known, the period map gives an isomorphism (of quasi-projective varieties): 
$$\calF_2\cong \calD/\Gamma,$$
where $\calD$ is a $19$-dimensional Type IV domain, and $\Gamma$ an arithmetic group. As discussed in \S\ref{compactdg},  $\calD/\Gamma$ has several compactifications, but neither of them is known to be geometric. Shah \cite{shah} has constructed a different model  (compact and with weak geometric meaning) for $\calF_2$ by using GIT. Specifically,  a degree $2$ $K3$ surface is a double cover of $\bP^2$ branched along a plane sextic. Thus, the GIT quotient $\overline \calM$ for plane sextics is a birational model for $\calF_2$. Not all degree $2$ $K3$ surfaces are double covers of $\bP^2$; the so called {\it unigonal} $K3$s (which form a divisor in $\calF_2$) correspond in $\overline \calM$ to the triple conic. More precisely, the triple conic is semistable with closed orbit; we denote by $\omega\in \overline \calM$ the corresponding point in the GIT quotient. Shah has shown that the exceptional divisor of a blow-up $\widehat\calM$ of $\omega\in \overline \calM$ (in modern language, the Kirwan blow-up of $\omega$) corresponds to the unigonal degree $2$ $K3$s.   It follows that the GIT space $\widehat\calM$ is a compactification of $\calF_2$ with  explicitly described boundary. 

From the Hodge theoretic perspective, $\calF_2\cong \calD/\Gamma$ has a canonical compactification, the SBB compactification $(\calD/\Gamma)^*$. In this situation the boundary of  $(\calD/\Gamma)^*$ consists of $4$ modular curves meeting in a point. On the other hand, the boundary of $\calF_2$ in the GIT compactification $\widehat\calM$ also has $4$ boundary components, but some of dimension $1$ and some of dimension $2$ (see figure \ref{fig4}). In other words, $\widehat{\calM}$ and $(\calD/\Gamma)^*$ almost agree everywhere. The explanation for this fact is given by the general framework of Looijenga discussed above. Specifically, the following holds:
\begin{figure}[htb!]
\includegraphics[scale=0.57]{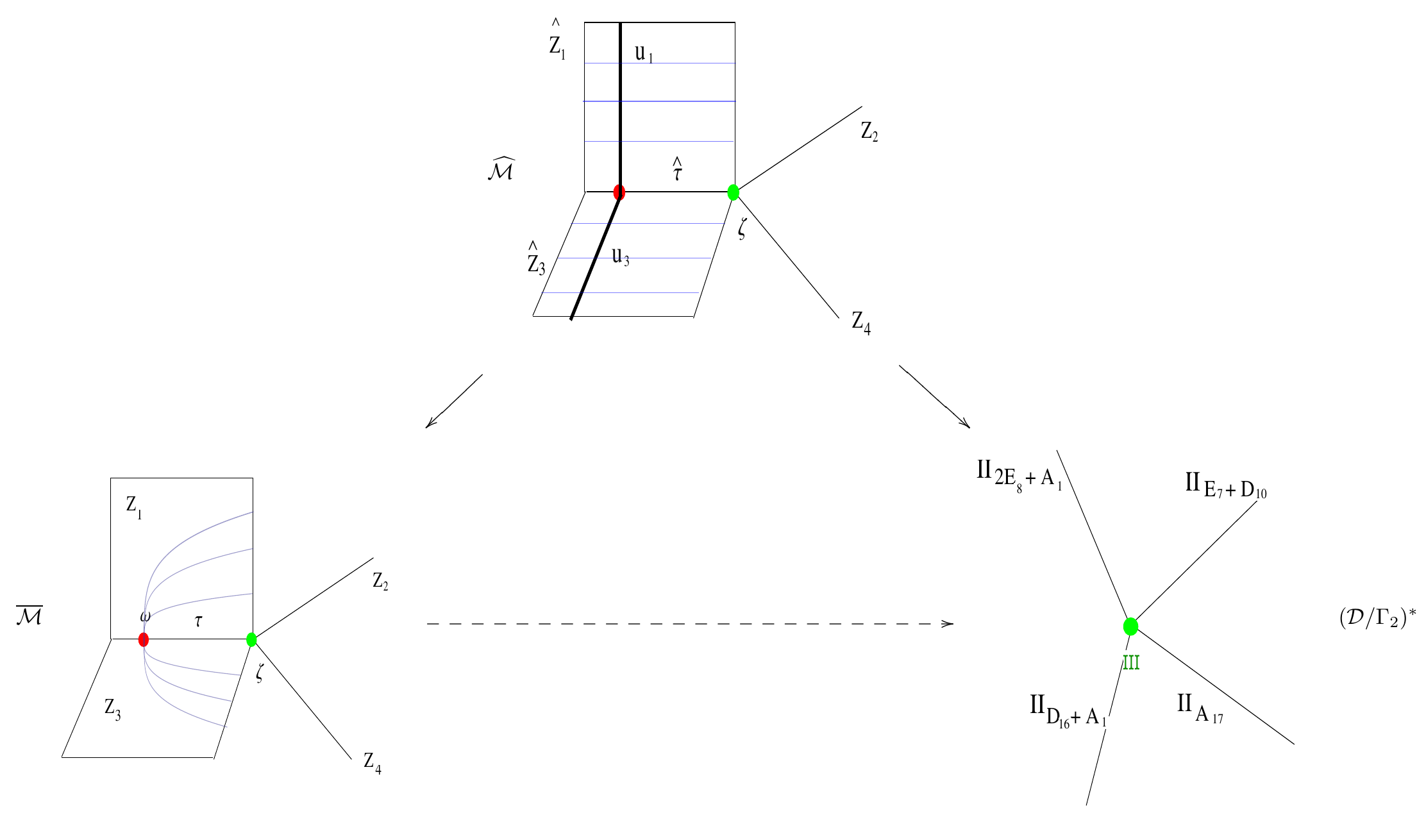} 
\caption{The boundary strata in the GIT and Baily-Borel compactifications for degree $2$ $K3$ surfaces}\label{fig4}
\end{figure}

\begin{theorem}[Looijenga, Shah] \label{thmlooijenga}
The open embeddings $\calF_2\subset \widehat{\calM}$ and $\calF_2\subset (\calD/\Gamma)^*$ extend to a diagram (with regular maps):
$$
\xymatrix{
&\widehat{\calM}\ar@{->}[ld]\ar@{->}[rd]&\\ 
\overline{\calM}\ar@{-->}[rr]&&(\calD/\Gamma)^*&
}
$$
such that
\begin{itemize}
\item[i)] $\widehat \calM\to \overline \calM$ is the partial Kirwan blow-up of $\omega\in \overline\calM$;
\item[ii)] $\widehat \calM\to (\calD/\Gamma)^*$ is the Looijenga modification of the Baily-Borel compactification associated to the hyperplane arrangement $\calH_\infty$ which corresponds to the divisor of unigonal $K3$s; more intrinsically, it is a small modification of $(\calD/\Gamma)^*$ such that the closure of Heegner divisor $\calH_\infty/\Gamma$ becomes $\bQ$-Cartier. 
\item[iii)] The exceptional divisor of $\widehat \calM\to \overline \calM$ maps to the unigonal divisor. 
\item[iv)] The boundary components are mapped as in Figure \ref{fig4}. 
\end{itemize}
\end{theorem}

\begin{remark}
We note that the birational map $\overline \calM\dashrightarrow (\calD/\Gamma)^*$ is everywhere defined except the point $\omega$. After the Kirwan blow-up $\widehat \calM\to \overline \calM$ of $\omega$, one gets a regular map. The explanation for this fact is that, except for the triple conic, all the degenerations occurring in the GIT quotient $\overline \calM$ are ``mild''. This was first observed by Shah \cite{shahinsignificant, shah}. In modern language, this is related to the theory of du Bois singularities (see \S\ref{sdubois}). 
\end{remark}

\section{The KSBA Approach to Moduli Spaces}\label{sectksba}
The previous two approaches to constructing moduli spaces are global in nature. Essentially, these approaches give a quasi-projective variety $\calM$ parameterizing the generic objects of a moduli problem, and  a compactification $\overline \calM$ with little control on the modular meaning of the boundary. The KSBA approach reviewed below takes the opposite point of view. It starts by identifying the correct boundary objects, and then it constructs a moduli space by gluing the local deformations of ``stable'' objects in a stack. It is of significant interest to reconcile this approach with the other two approaches (thus gaining modularity at the boundary, as well as nice global properties). This topic is still in early stages of development. We briefly discuss below the du Bois singularities which link to the Hodge theoretic approach, and K-stability which links to GIT.  

The original KSBA approach was introduced in Koll\'ar--Shepherd-Barron \cite{ksb} and Alexeev \cite{alexeev}. A recent survey of the KSBA approach is Koll\'ar \cite{kollar}. A two volume monograph (see \cite{ko13}) handling the more delicate aspects of the KSBA approach will appear soon. A detailed example (including a discussion of the state-of-art theoretical background) is discussed in Aleexev's lectures \cite{alexeevbarcelona}.
\subsection{The KSBA approach}
The KSBA approach is a natural generalization of the Deligne--Mumford construction  \cite{dm} of the proper stack of stable curves $\overline\calM_g$. We review this approach below focusing on the main ideas of the construction and ignoring various technical issues. We caution the reader that some of these technical are quite subtle in higher dimensions and should consult \cite{kollar,ko13} for precise statements. 

\subsubsection{The original construction of $\overline \calM_g$ does not involve GIT.} The abstract point of view to moduli is to define a {\it moduli functor}, i.e. to define a class of objects of interest and allowable families of objects (e.g. flat proper morphisms $\calX/S$). Under mild hypotheses, there is an associated algebraic moduli stack  $\calM$, which is roughly obtained by gluing the (allowable) deformation spaces (= local patches of $\calM$) of the objects under consideration. From this point of view, the main issue is to select a class of objects so that the moduli stack is proper and separated. This is a delicate balance, as it is easy to achieve either of those two conditions, e.g. smooth varieties of general type (with ample $K_X$) vs. GIT semistable varieties (say with respect to the embedding given by $nK_X$), but not both at the same time. By the valuative criteria, defining a proper and separated moduli stack is essentially equivalent to asking that a $1$-parameter family $\calX^*/\Delta^*$ over the punctured disk has a unique limit with respect to the given moduli functor.

In the case of curves, it is well known that the correct class of objects to define a proper moduli stack $\overline\calM_g$ is that of {\it stable curves} $C$:
\begin{itemize}
\item[(S)]  $C$ is a curve with only nodal singularities;
\item[(A)] every rational component contains at least 3 special points, or equivalently $|\Aut(C)|<\infty$,  or equivalently $\omega_C$ is ample. 
\end{itemize}
This follows from considering the following three steps:
\begin{itemize}
\item[Step 0:] {\it Semi-stable reduction} - after a finite base change we can assume $\calC/\Delta$ is semistable family (i.e. $\calC$ is smooth with reduced normal crossing central fiber). This is a nice topological and geometric model, which is the starting point for the next steps. Of course, the issue here 	is that $\calC$ is far from unique (e.g. due to blow-ups). 
\item[Step 1:] {\it Relative minimal model} - $\calC$ is a surface, and thus the minimal model can be easily obtained by contracting the $(-1)$-curves in the fibers. Here the total space $\calC$ is still smooth, but the moduli functor is not yet separated. 
\item[Step 2:] {\it Relative canonical model} - this is obtained by contracting the $(-2)$-curves orthogonal to $\omega_{\calC/\Delta}$. Here, we obtained separateness, but we need to allow mild singularities (rational double points) for the total space.  
\end{itemize}
\begin{remark}
The main facts leading to nodal curves as the correct class of degeneration for curves are: (a) the relative canonical model agrees with the log canonical model of the pair $(C, C_0)$, and then (b) by adjunction, it follows that $C_0$ has semi-log canonical singularities and hence $C_0$ is nodal.
\end{remark}

In general, Step 0 is true in high generality (cf. \cite{mumfordtai}). 
Steps 1 and 2 are then achieved (under relatively mild assumptions) by using the full strength of the Minimal Model Program (MMP; \cite{bchm}). We emphasize that for Step 2 it is essential to have varieties of general type (or  log general type). For instance, we note that for $K3$ surfaces there exists a very good answer to steps 0 and 1:
\begin{theorem}[Kulikov--Persson--Pinkham/Shepherd-Barron \cite{k3book}] Let $(\calX,\calL)/S$ be a $1$-parameter degeneration of polarized $K3$ surfaces. After a finite base change and birational transformations, one can assume the following:
\begin{itemize}
\item[i)] $\calX$ is a semistable family;
\item[ii)] $K_{\calX/B}\equiv 0$;
\item[iii)] $\calL$ is nef. 
\end{itemize}
\end{theorem}
\noindent However, there is no good analogue to Step 2. Namely, one can define $\overline \calX=\Proj R(\pi_*\calL)$, which is quite similar to a relative log canonical model, but even so one does not obtain a separated moduli stack. 
If one considers instead of the line bundle $\calL$ a divisor $\calH$ (with $\calL=\calO_\calX(\calH)$) the moduli problem becomes a problem about surfaces of log general type and it has a good solution (see \cite{lowdeg}).

\subsubsection{Due to MMP, the Deligne--Mumford approach ``easily'' generalizes in higher dimensions.}
As discussed above, for varieties of (log) general type the Steps 0-2 from the construction of $\overline \calM_g$ go through also in higher dimensions. Analogous to case of curves, the correct limits of smooth varieties of general type are the KSBA {\it stable varieties} $X$ that satisfy two conditions analogue to the stable curve conditions:
\begin{itemize}
\item[(S)] $X$ has semi-log-canonical singularities (or more precisely $X$ is an slc  variety, see \cite[Def. 3.1]{kollar});
\item[(A)] $\omega_X$ is ample. 
\end{itemize}
Note that the similarity to curves is somewhat deceiving. There are several subtle differences that make the results in higher dimension to be much harder than those for curves. For instance, even to define the conditions (S) and (A), one needs to assume: $K_X$ is $\bQ$-Cartier (N.B. the nodal curves are Gorenstein; so the issue does not arise).  Similarly, while all stable curves are smoothable, this is not true for stable varieties in higher dimensions. This leads to various issues: If one restricts to the smoothable stable varieties, there is no natural scheme structure at the boundary. On the other hand, if one considers all KSBA stable varieties, the application of  MMP as sketched above can fail on the components of the moduli that parameterized non-smoothable slc varieties. A different method to handle these components is needed (see \cite[\S5.2]{kollar}).

Another set of issues in higher dimensions is the correct definition of the moduli functor (i.e. admissible families). While for curves, the natural requirement of flat and proper morphism suffices, in higher dimensions one has to impose additional conditions (such as the total space of a family has to be $\bQ$-Gorenstein). We refer to \cite[(4.3),(4.4)]{kollar} for a discussion of the moduli functor. Some of the more subtle issues were only recently settled. We refer the interested reader to \cite{kollar} and the upcoming companion monograph. 

The upshot of this discussion is that, due to deep results in birational geometry and substantial additional work, there exists a proper and separated moduli stack $\calM_h$ for varieties of (log) general type (with Hilbert polynomial $h$)  which generalizes $\overline{\calM}_g$.  Also, since a KSBA stable variety has no infinitesimal automorphisms (e.g. \cite[Lemma 2.5]{zsolt}), $\calM_h$ is a Deligne-Mumford stack. Consequently $\calM_h$ has associated an algebraic space as coarse moduli space (\cite{km}). Unfortunately, as discussed below the analogy to  $\overline{\calM}_g$ stops here. 

\subsubsection{The only property of $\overline M_g$ that generalizes in higher dimensions is the projectivity.} 
The moduli stack  $\overline{\calM}_g$ is a proper smooth Deligne--Mumford stack, with a coarse projective moduli stack $\overline M_g$ (which, as previously discussed, can be constructed via GIT).   Additionally, the boundary of $\overline{\calM}_g$ is a divisor with normal crossings. Unfortunately, no general smoothness result can hold for the moduli of varieties of general type. The reason for this is that the deformation spaces even for smooth varieties can behaved arbitrarily badly (i.e. on the versal deformation spaces one can encounter essentially any singularity that can be defined over $\bZ$). This is the content of Vakil's results:
\begin{theorem}[Vakil \cite{vakilmurphy}] 
The following moduli spaces satisfy Murphy's law.
\begin{itemize}
\item[M2a.] the versal deformation spaces of smooth n-folds (with very ample
canonical bundle), $n\ge 2$.
\item[M3.] the Hilbert scheme of nonsingular surfaces in $\bP^5$, and the Hilbert scheme of surfaces in $\bP^4$. 
\end{itemize}
\end{theorem}
The first result says that arbitrarily bad singularities can occur even for moduli of smooth surfaces of general type.  The second result is relevant for the GIT constructions. Of course, as previously discussed, the GIT compactifications  will tend to have even worst singularities than the KSBA compactifications. 

One deep positive result is Koll\'ar projectivity \cite{kollarp} of the coarse moduli space (see also \cite{fujino}): 
\begin{theorem}[Koll\'ar, Fujino]\label{projmodt}
The moduli functor of  stable $\calM_h$ varieties with Hilbert function $h$ is coarsely represented by a projective algebraic scheme.
\end{theorem}

\begin{remark}\label{projmod}
The methods of Koll\'ar use essentially that $\calM_h$ is proper.  It is not possible to prove (by those methods) the quasi-projectivity of the smooth locus without first compactifying the moduli space. On the other hand, as already noted, Gieker \cite{gieseker}  proved (using GIT) the quasi-projectivity of the moduli space of smooth surfaces of general type. A similar higher-dimensional result (but via significantly different GIT methods) was proved by Viehweg (\cite{viehweg}). However, even in these quasi-projectivity results, there is a compactification (i.e. a GIT compactification) in the picture.   In other words, we see here a manifestation of the well known principle in algebraic geometry: {\it even if we only care about the smooth locus, it is important to have a compactification.}
\end{remark}

\subsubsection{It is quite difficult to understand the KSBA compactification. In fact, few concrete examples are known. The GIT approach might help.} 

The bad behavior of the moduli of varieties of general type discussed in the previous section has to do more with the pathology of the deformation spaces rather than the moduli itself. One might hope that for nice classes of varieties there will be a good moduli space with applications similar to those of $\overline M_g$. Unfortunately, very few  examples of KSBA compacifications for moduli spaces are known. Namely, while a number of interesting  illustrations   of the KSBA compactification procedure were given  (e.g. \cite{hacking},  \cite{opstall}, \cite{tevelev},  \cite{pardini}),  all  these examples tend to be special in the sense that they are  related to curves or hyperplane arrangements, etc. It is of interest to give more ``generic'' examples (such as moduli of quintic surfaces, or other surfaces with small $p_g$ and $q=0$) of compactifications for moduli of surfaces of general type. 

The KSBA compactification is clearly the ``correct'' compactification for varieties of general type. However, its abstract definition makes it somewhat intractable. (What is special for $\overline \calM_g$ is that a stable genus $g$ is obtained from lower genus curves via a simple combinatorial gadget, the dual graph, and then one can proceed inductively. For surfaces and higher dimensions, the situation is exponentially more involved.) Part of the message of these notes is to advocate that it might be possible to arrive at a KSBA compactification via interpolation. For instance, for quintic surfaces, a GIT compactification is readily available. In general, 
for GIT the numerical criterion gives an algorithmic way of determining the semistable points. No algorithmic procedure is known for KSBA. Thus, GIT tends to be more accessible. To interpolate from the GIT to the KSBA compactification, one can apply the replacement lemma \ref{ssreplace} and reduce to the study of $1$-parameter families $\calX/\Delta$ with semistable (even with closed orbit) central fiber $X_0$ but not log canonical. The question then is to find the KSBA replacement for $X_0$ (for related work for curves see \cite{hstable}). 

Returning to the example of quintic surfaces, Patricio Gallardo \cite{patricio} and Julie Rana have shown that the KSBA replacement for quintic surfaces $X_0$ with a unique Dolgachev singularity (those are GIT stable) are of type $\widetilde X_0\cup S$, where $\widetilde X_0$ is the resolution of $X_0$ and $S$ is a $K3$ surface ($S$ is a ``$K3$ tail'', which depends on the degeneration). These type of surfaces give  divisors in the KSBA compactification for quintics. It is interesting to note that these type of example were encountered also in work of Friedman \cite{friedman5} and Shepherd-Barron who were concerned with the Hodge theoretic properties of degenerations of quintics. One wonders if the period map approach to moduli would not shed further light on the compactification problem for quintics by bringing in tools such as representation theory.

\begin{remark}\label{remlogcan}
A simple connection between GIT and KSBA was observed by Hacking \cite{hacking} and Kim-Lee \cite{kimlee}. Namely, for hypersurfaces $V\subset \bP^n$ of degree $d$: {\it if $(\bP^n, \frac{n+1}{d} V)$ is a log canonical pair, then $V$ is GIT semistable.} Similarly, we have {\it if $(\bP^n, (\frac{n+1}{d}+\epsilon) V)$ (for $0<\epsilon\ll 1$) is a log canonical pair, then $V$ is GIT stable}. The converse is not true, an example being given by the triple conic viewed as a plane sextic. The reason for this is that the log canonical threshold and the numerical function $\mu(x,\lambda)$ are computed in the same way by minimizing a certain quantity over all choices of coordinates. The difference is that for GIT only linear changes of coordinates are allowed vs. analytic changes of coordinates in the KSBA case. In conclusion,  GIT allows worst singularities (which are not detected by linear coordinates) for the objects in the boundary and consequently more collapsing (resulting in loss of geometric information on the degenerations).
\end{remark}

\subsection{Slc singularities are du Bois (connection to Hodge theory).}\label{sdubois} 
From a topological/Hodge theoretic point of view, the good models for studying degenerations are the semistable models $\calX/\Delta$. In this situation, the Clemens--Schmid exact sequence (e.g. \cite{clemensschmid}):
\begin{equation}\label{eqcs}
\dots \to H^n(X_0)\to H^n_{\lim}\xrightarrow{N} H^n_{\lim} \to H_n(X_0)\to \dots
\end{equation}
(where $N=\log T$ is the monodromy, and $H^n_{\lim}$ is the limiting mixed Hodge structure)  says that the limiting Hodge structure is essentially determined by the mixed Hodge structure of the central fiber $X_0$. The KSBA limits are obtained from $\calX/\Delta$ by considering the relative log canonical model $\calX^c/\Delta$. It is natural to ask what is the connection between the KSBA limits and the Hodge theoretic limits. In particular, it is interesting to understand the connection between the MHS on $X_0^c$ and $H^n_{\lim}$. If these questions have satisfactory answers, then, in principle, by using the structure on the Hodge theoretic side (see Section \ref{sectht}) one gets information of the KSBA compactification. Unfortunately, since $X_0^c$ and $\calX^c$ are quite singular, the situation is not completely understood. We described below a step towards a comparison between KSBA and Hodge limits, but much remains to be done.

Since the KSBA stable curves are nodal and $\overline\calM_g$ is a smooth normal crossing compactification, the KSBA and Hodge theoretic models are closely related in dimension $1$. Namely, as discussed in \S\ref{compactdg}, 
 there is an extended Torelli morphism
$$\overline{\calP}:\overline{M}_g\to \overline{\calA}_g^{Vor}$$
 which corresponds to class $\lambda$ on $\overline{M}_g$ (and thus $\lambda$ is semi-ample). For surfaces, Mumford and Shah \cite{shahinsignificant} noted that a certain class of singularities (i.e. Gorenstein slc singularities in modern language) are {\it cohomologically insignificant singularities}. Steenbrink \cite{steenbrinkinsignificant} then interpreted this statement as saying that (in dimension $2$) the Gorenstein slc singularities are {\it du Bois singularities}. The du Bois singularities are a class of singularities that behave well from a Hodge theoretic point of view. Roughly speaking, if the central fiber $X_0^c$ has du Bois singularities, then the natural comparison map
 \begin{equation}\Gr_n^W\mathrm{sp}_n:\Gr_n^WH^n(X_0^c)\to \Gr_n^WH^n_{\lim}
 \end{equation}
  is an isomorphism. Shah \cite{shah} has used this fact to study the moduli of degree $2$ $K3$ surfaces (see the discussion of \S\ref{githodge}). 

It turns out that the results for surfaces (slc $\Longrightarrow$ du Bois) generalize well. Specifically, Koll\'ar and Kov\'acs have obtained the following general result:
\begin{theorem}[{Koll\'ar--Kov\'acs \cite[Thm. 1.4]{kk}, \cite[Cor. 6.32]{ko13}}]
Let $(X,\Delta)$ be an slc pair. Then $X$ is du Bois. 
\end{theorem}
 
 As noted, one application of the previous theorem is to understand images of period maps (see \S\ref{githodge}).  Conversely, using Hodge theory, one obtains results on the KSBA moduli. For instance, the following results, which settles several technical issues in the KSBA construction, is obtained by using du Bois singularities. 
\begin{theorem}[Koll\'ar, Kov\'acs, cf. {\cite[p. 263]{ko13}}]
Let $f:\calX\to S$ be a proper and flat morphism with slc fibers over closed points; $S$ connected. 
\begin{itemize}
\item[(1)] Let $\calL$ an $f$-semi-ample line bundle on $\calX$ (e.g. $\calL=\calO_\calX$). Then, for all $i$, 
\begin{itemize}
\item[(a)] $R^if_*(\calL^{-1})$ is locally free and compatible with base change and
\item[(b)] $h^{i}(X_s,L_s^{-1})$ is independent of $s\in S$.
\end{itemize}
\item[(2)] If one fiber of $f$ is Cohen-Macaulay then all fibers of $f$ are Cohen-Macaulay.
\item[(3)] $\omega_{\calX/S}$ exists and is compatible with base change. Furthermore, for all $i$, 
\begin{itemize}
\item[(a)] $R^{i}f_*\omega_{\calX/S}$ is locally free and compatible with base change and
\item[(b)] $h^i(X_s,\omega_{X_s})$ is independent of $s\in S$.  
\end{itemize}
\end{itemize}
\end{theorem}

\subsection{Asymptotic stability, K-stability, and KSBA}\label{kstability}
We close our survey with a brief discussion of the connection between (appropriate) GIT stability and KSBA stability. For further details and a more technical discussion, we recommend the survey \cite{odakasurvey}.

\subsubsection{Asymptotic stability and KSBA}
The following theorem can be regarded as a comparison theorem between the GIT and KSBA approach. It says that if a canonically polarized variety $X_0$ is asymptotically semistable (even in a weak sense) then it is also KSBA stable.

\begin{theorem}[Wang-Xu \cite{wang}]
Let $\calX/S$ be a KSBA-stable family over a smooth curve $S$, and $o\in S$ a special point.  Assume that the fibers $X_s$ for $s\in S^{\circ}=S\setminus \{ o\}$ are asymptotically (Chow) stable. Let $r$ be such that $rK_\calX$ is Cartier. Then for any flat family $(\calX',\calL)/S$ of polarized varieties which agrees with $(\calX,\omega^{[r]}_\calX)$ over $S^{\circ}$ and which has the property that the special fiber $(X_o,L_o^k)$ is Chow semistable for infinitely many $k>0$, we have $(\calX', \calL)\cong (\calX,\omega^{[r]}_{\calX})$ (i.e. agrees everywhere).
\end{theorem}

Unfortunately, the converse can not hold. It was observed by Mumford (\cite[Prop. 3.4]{mumford}, \cite[Prop. 3.1]{shahinv}) that 
\begin{proposition}[Mumford]
An asymptotically semistable surface has singularities of multiplicity at most  $6$.
\end{proposition} At the same time, there are KSBA stable surfaces with singularities of higher multiplicity. Together with the theorem above, this implies a negative result:
\begin{corollary}
Assume $X_0$ is a KSBA stable surface which is not asymptotically semistable. Then $X_0$ has no asymptotic semistable replacement.  
\end{corollary}

In other words, the asymptotic GIT quotients will not stabilize (as the embedding gets higher and higher); for some concrete examples of this see \cite{wang} (also \cite{vedova}).

\begin{remark} In general, it is very hard to analyze the asymptotic stability. The best positive result so far is a hard  result of Gieseker which says that 
a smooth surface of general type is asymptotically stable. More precisely, the following holds:

\smallskip

\noindent{(Gieseker \cite{gieseker})} {\it Let $X$ be a surface of general type. For any sufficiently large $n$, there is an $m$ so that the $m^{th}$ Hilbert point of the $n$ canonical image of $X$ is stable. Furthermore, $m$ and $n$ depend only on $K^2$ and $\chi(\calO_X)$.}
\end{remark}

\subsubsection{K-Stability}
In recent years, motivated by the work of Donaldson, Tian, and Yau on the existence of special metrics, a new notion of stability, {\it K-stability}, has emerged (e.g. \cite{donaldson2}). Namely, conjecturally, {\it a polarized K\"ahler manifold $(X,L)$ admits a constant scalar curvature metric with class $c_1(L)$ iff $(X,L)$ is $K$-polystable.} This is an area of very active research with many important recent results, but this goes beyond the purpose of this survey. We only touch here on the connection between K-stability and algebraic geometry. This is mostly based on work of Odaka (see \cite{odakasurvey} for a related survey).

 The K-stability can be viewed as a refined notion of asymptotic stability. 
In fact, the set-up for K-stability is quite similar to that of the numerical criterion from GIT. Namely, one considers {\it test configurations}, i.e. $1$-parameter families $\calX/\bA^1$ equivariant with respect to a $\bC^*$ action (with $\bC^*$ acting in the standard way on $\bA^1$, and $X_t\cong X_{t'}$ for $t,t'\neq 0$). To such a test configuration, one associates a numerical invariant, the Donaldson--Futaki invariant $DF(\calX)$ (analogue to $\mu(x,\lambda)$). Then, as in the GIT numerical criterion \ref{numcriterion}, the K-(semi)stablity is equivalent to the positivity/non-negativity of the Donaldson--Futaki invariant  for all non-trivial test configurations. 

It turns out that K-stability is essentially equivalent to KSBA stability. Specifically, the following hold:

\begin{theorem}[{Odaka \cite[Thm. 1.2]{odaka}}] Let $X$ be a projective scheme satisfying (*) (see \cite[Def. 1.1]{odaka}) and $L$ be an ample line bundle on $X$. Then, if $(X,L)$ is K-semistable, $X$ has slc singularities.
\end{theorem}
\noindent(The condition (*) encodes the pre-conditions on $X$ to even make sense to say that $X$ has slc singularities.)
\begin{theorem}[{Odaka \cite[Thm. 1.5]{odaka}}] The following hold:
\begin{itemize}
\item[i)] A slc polarized variety $(X,L)$ with numerically trivial canonical divisor $K_X$ is K-semistable.
\item[ii)] A slc canonically polarized variety $(X,K_X)$ is K-stable. 
\end{itemize}
\end{theorem}

The main drawback of K-stability is that it is not known how it can be used to construct a moduli space. For instance, it is not known if it is an open condition. Regardless, it can be viewed as an interpolation between two algebro-geometric stability conditions: GIT and KSBA stability. One wonders if 
it might be possible to set-up a slight modification of the asymptotic stability to obtain K-stability (and thus KSBA stability). An example (in spirit) going in this direction is \cite{odakasun}, but much remains to be done.

\bibliography{moduli}
\noindent Stony Brook University, Department of Mathematics,  Stony Brook, NY 11794.
\\ {\it E-mail address}: \tt{rlaza@math.sunysb.edu}
\end{document}